\documentclass[12pt]{amsart} 
\newfont{\sheaf}{eusm10 scaled\magstep1}

\newcommand{\ra}{\ensuremath{\rightarrow}}

\def\eea{\end{eqnarray*}}
\def\bea{\begin{eqnarray*}}

\newcommand{\Proof}{{\it Proof. }} 
\newcommand{\QED}{{\hfill $Q.E.D.$}} 
 
\newtheorem{teo}{Theorem}[section] 
 
\newtheorem{lem}[teo]{Lemma}

\newtheorem{oss}[teo]{Remark}

\newcommand{\C}{\ensuremath{\mathbb{C}}}
\newcommand{\R}{\ensuremath{\mathbb{R}}} 
\newcommand{\Z}{\ensuremath{\mathbb{Z}}}

\newcommand{\NNN}{\ensuremath{\mathcal{N}}}
\newcommand{\F}{\ensuremath{\mathbb{F}}}

\newcommand{\I}{\ensuremath{\mathcal{I}}}

\newcommand{\N}{\ensuremath{\mathbb{N}}} 
\newcommand{\hol}{\ensuremath{\mathcal{O}}}
 
\newcommand{\PP}{\ensuremath{\mathbb{P}}}

\newcommand{\SSS}{\ensuremath{\mathcal{S}}}

\begin{document}

\title{The 3-cuspidal quartic and  braid monodromy of   degree 4 coverings.}

\author{Fabrizio Catanese - Bronislaw Wajnryb\\
}
 \footnote{
The research of the  author was performed in the realm  of the 
 SCHWERPUNKT "Globale Methode in der komplexen Geometrie",
and of the EAGER EEC Project. \\

}
 \\

\date{October 29, 2004}
\maketitle

\begin{abstract}
Motivated by the study of the differential and symplectic topology
of $(\Z/2)^2$- Galois covers of $\PP^1\times \PP^1$, we determine
the local braid monodromy of natural deformations of
smooth $(\Z/2)^2$- Galois covers of surfaces at the points where
the branch curve has a nodal singularity.

The study of the local deformed branch curves is solved via some interesting
 geometry of projectively unique objects: plane quartics with $3$  cusps, 
which are the plane sections of the  quartic surface having the twisted 
cubic as a cuspidal curve.
\end{abstract}

\section{Introduction}

 This article is a continuation of a preceding one (\cite{c-w}), which was devoted to
the proof that the so called $(a,b,c)$-surfaces (where we take $a,b,c \in \N $
 with $b$ and $ a+c$ fixed) provide
examples of simply connected algebraic surfaces which are
diffeomorphic but not deformation equivalent.

The  $(a,b,c)$-surfaces are coverings of $\PP^1\times \PP^1$ of degree $4$
and are defined by 2
equations
\begin{eqnarray}  z^2 &=&
  f(x,y)\\
  w^2 &=&
  g(x,y) \nonumber ,\end{eqnarray}

where f and g are bihomogeneous polynomials ,  belonging to
  respective vector spaces of sections of line bundles:
  $ f \in H^0({\PP^1\times \PP^1}, {\hol}_{\PP^1\times \PP^1}(2a,2b)) $ and

$ g \in H^0({\PP^1\times \PP^1}, {\hol}_{\PP^1\times \PP^1}(2c,2b)). $
\bigskip

A question which was left open in \cite{c-w} was the symplectic equivalence
of the above $(a,b,c)$-surfaces. To this purpose, and for more general purposes,
it is important to determine the braid monodromy factorization of the branch curve 
corresponding
to a symplectic deformation of the 4-1 covering $S \ra \PP^1 \times \PP^1$
possessed by an $(a,b,c)$-surface $S$
(note that in \cite{c-w} one key result was the determination of the 
mapping class group monodromy factorization, which is a homomorphic image
of the braid monodromy factorization).

In this paper we approach the first step, namely, we determine the local 
braid monodromy factorization of 4-1 coverings which are deformations
of bidouble covers ($(\Z/2)^2$- Galois covers).

The discriminant picture that we get is somehow unexpected, in that instead of the
usual swallowtail surface we obtain a rational quartic surface with a twisted
cubic as cuspidal curve. We show in the last section, namely in Theorem \ref{edge}, 
that such surface is projectively unique, being the tangential developable of
the twisted cubic, or equivalently, the dual surface of the
twisted cubic curve in $\PP^3$, or the discriminant of the general equation of degree $3$.

As hinted at in the first section, the picture is
not completely unexpected, especially the fact that the quartic is a discriminant surface
for the general equation of degree $3$, since actually 
 Galois theory teaches us that deformed Galois
covers of degree 4 are exactly the trick to relate the 
solvability of the general equation of degree 4 to the solvability of
the general equation of degree 3.
It follows that our perturbed local discriminant curve of
 $S \ra \PP^1 \times \PP^1$ is a plane  quartic curve $\Delta$ with three cusps,
and bitangent to the line at $\infty$ in two real points.

Viewing the curve as a small perturbation of a pair of real lines counted with
multiplicity two made the determination of the braid monodromy 
of this affine curve almost impossible.

The trick which solved the problem is the following well known observation:
a three cuspidal quartic over an algebraically closed field is projectively
unique, since it is the dual curve of a nodal plane cubic curve.

We can then change the real picture and take a nodal cubic with an isolated
double point, but with three real flexes: its dual curve, once we take as line 
at infinity the dual line of the nodal point, will be a quartic $C$ with three
real cusps, and bitangent at the line at infinity in two imaginary points. 

For $C$ the  points with a real abscissa $x$ which are interesting
for the determination of the braid monodromy have now an ordinate $y$ which is 
either real or imaginary, and it is quite easy to calculate then the 
braid monodromy factorization. From this one, since $\Delta$ is complex affine
equivalent to $C$, we deduce the braid monodromy factorization for $\Delta$.

Lack of time prevents us to analyse the question whether one can similarly determine the 
local braid monodromy factorization for a deformed abelian cover. This would also
be a very useful result in the study of the differential and symplectic
topology of huge classes of algebraic surfaces.


\section{The discriminant of a deformed bidouble cover. }

Consider a ring $A$ of characteristic $p \neq 2$ and a so called simple bidouble cover of
$A$, 
i.e., a $(\Z/2)^2$- Galois ring extension
$ A \subset B'$  where $B'$ is the quotient ring of $A [z,w]$ given, for 
some choice of $u,v \in A$, by

\begin{eqnarray}  z^2 &=&
 v\\
 w^2 &=& u \nonumber .\end{eqnarray}

A {\bf deformed bidouble cover} is a finite ring extension 
$ A \subset B$ given, for $u,v,a,b \in A$, by

\begin{eqnarray}  z^2 &=&
 v + a w \\
 w^2 &=& u + bz\nonumber .\end{eqnarray}

Observe that, if $a$ is invertible, then  $ w = a^{-1} (z^2 - v)$, and
we get the quartic equation  $  a^{-2} (z^2 - v)^2 - bz - u = 0$, equivalent
to $   z^4 - 2 z^2 v  -  a^2 b z  + (v^2 -  a^2 u) = 0 $ and that,
since $p \neq 2$, every quartic equation can be reduced, by a translation
(Tschirnhausen transformation), to the above equation, for a suitable
choice of $v, b, u$.

This  standard trick, which allows to deform a bidouble
Galois extension to the general quartic equation is very important
in Galois theory. Because, in the $z,w$ plane we have the pencil of conics
generated by the above two parabolae, and the determinant function on
the parameter of the pencil provides an equation of degree three,
whose Galois group corresponds to the image of the Galois group of the quartic equation
under the surjection $\SSS_4 \ra \SSS_3$ with kernel
$(\Z/2)^2$.

 We are however interested in a finer geometric question, we do not only look at
 algebraic extensions of function fields, indeed we look more closely
 at finite coverings of smooth algebraic varieties.

 The concepts of bidouble cover, and deformed bidouble cover have been introduced, in
the global case of coverings of smooth algebraic varieties, in \cite{cat1}, the latter
under the name of {\bf natural deformations} of bidouble covers. 
We refer the reader for details to \cite{cat1}, and as well  to \cite{c-w}
for the applications we have in mind.

We observe that $B$ is a rank 4 free $A$-module, with basis $1, z,w,zw$, and
that to each nontrivial basis element corresponds the respective multiplication 
matrix

$$M_z=
\begin{pmatrix}
  0 &v& 0 & a u \\
  1  &0&0 & ab  \\
   0 & a &0 & v  \\
   0 & 0& 1 &0\\

\end{pmatrix}, 
M_w=
\begin{pmatrix}
  0 &0& u & bv \\
  0 &0&b & u  \\
   1 & 0 &0 & ab  \\
   0 & 1& 0 &0\\

\end{pmatrix},
M_{zw}=
\begin{pmatrix}
  0 &au& bv & uv \\
  0  & ab & u & bv  \\
   0 & v & ab  & au  \\
   1 & 0& 0 &ab\\

\end{pmatrix}.$$

Consider now the different $R$ (ramification Cartier divisor) of the ring extension, 
$$ R = det \begin{pmatrix}
  2z &-a \\
  -b &2w \\
\end{pmatrix} = 4 zw - ab.$$

 We can then find the discriminant $\Delta$ as the norm of $R$, thus 

$$ \Delta = 4^4 det (M_{zw} - \frac{ab}{4} Id ) = 4^4 det \begin{pmatrix}
- ab/4 &au& bv & uv \\
  0  & 3ab/4 & u & bv  \\
   0 & v & 3ab/4  & au  \\
   1 & 0& 0 & 3ab/4 \\

\end{pmatrix} ,$$and
$$ - \frac{1}{16^2}\Delta =  - u^2 v^2 - \frac{9}{8} uv (ab)^2 + b^2 v^3 + a^2 u^3 + \frac{27}{16^2}
a^4 b^4  : = P(u,v,a^2,b^2).$$ 

We see immediately that, setting $\alpha := a^2, \beta := b^2$,
  $P(u,v,\alpha, \beta)$ is  homogeneous of degree $4$, and symmetric for
the involution $ (u, \alpha) \leftrightarrow (v, \beta ) $,
whose fixed point locus is not contained in $\{P =0\}   $ 
(this symmetry is forced by the symmetry exchanging 
$a$ with $b$, $w$ with $z$, $u$ with $v$).

\begin{teo}\label{disc}
The quartic hypersurface $P \subset \PP^3$ defined by $P(u,v,\alpha, \beta) =0$
is irreducible and has as singular locus a twisted cubic curve $\Gamma$, which is
a cuspidal curve for $P$. In particular $P$ is projectively unique,
being the dual surface of the twisted cubic. $P$ is also the tangential developable of the
twisted cubic, and the discriminant surface of the space $\PP^3$ of polynomials
of degree $3$ on $\PP^1$.
\end{teo}

\begin{oss}
Chapter V of \cite{supraz} is devoted to more general surfaces of degree 4
in $\PP^3$ which have a twisted cubic as double curve. 
\end{oss}

\Proof

We calculate for later use 
$$\partial P / \partial u = -2 u v^2 - \frac{9}{8} (ab)^2 v + 3 a^2 u^2 , \
\partial P / \partial v = -2 v u^2 - \frac{9}{8} (ab)^2 u + 3 b^2 v^2$$
and moreover 
$$\partial P / \partial 
\alpha = - \frac{9}{8}  u v (b)^2  + u^3  + \frac{54}{16^2}
a^2 b^4 , \
\partial P / \partial \beta =  - \frac{9}{8}  u v (a)^2 u + v^3
 + \frac{54}{16^2}
a^4 b^2,$$
whence in particular $$ u (\partial P / \partial u) - v (\partial P / \partial v) =
3 a^2 u^3 - 3 b^2 v^3,\alpha (\partial P / \partial 
\alpha) - \beta  (\partial P / \partial \beta) = a^2 u^3 -  b^2 v^3.$$

We conclude that the hypersurface $P$, i.e., $\{P =0\}   $, is irreducible, 
being reduced and being the image of the
quadric $ Q := \{  4 zw - ab = 0 \}$. We claim now that $P$ is singular along a 
twisted cubic $\Gamma$, which is of cuspidal type.

In view of the irreducibility of $P$, we will then conclude that 
there is no other singular curve on $P$.

In order to do this, let us work with affine coordinates on $Q$
 setting $ a=1$, whence $ b = 4 zw $ and $z,w$ are affine coordinates.

In terms of these coordinates,  $\alpha = 1$, $\beta = 16 z^2 w^2$, $u = w^2 - 4 z^2 w$,
$ v = z^2 - w$.

We conclude that $ F : Q \ra P$ factors through $ (z,w) \ra (s := z^2, w)$
and $G(s,w) : = (16 s w^2, w^2 - 4 s w, s-w)$.

 We  calculate the derivative matrix of $G$,
$$ DG = \begin{pmatrix}
  16 w^2 &32 s w  \\
  - 4w  & 2 w - 4 s  \\
   1& -1  \\
\end{pmatrix} $$

which has rank equal to 1 exactly for $ w + 2s=0$. Observe moreover that in these
points the kernel of $DG$ is given by the tangent vector $\partial/ \partial s + 
\partial/ \partial w $.

An immediate calculation shows that the image curve $\Gamma$ is the twisted cubic
$(64 s^3, 12 s^2, 3 s) $ and one may verify that on $\Gamma$  all the four partial
derivatives of $P$ do vanish.

It follows that $\Gamma$ is the only singular curve of $P$ (an irreducible quartic curve 
has at most three singular points), and that it is a cuspidal curve, since 
if we intersect $P$ with a general plane, then we get a curve in the $(s,w)$
plane which, at an intersection point with $ w + 2s = 0$, is tangent to the kernel of $DG$,
but maps with local degree one (since , for instance, the line
$ w = s + c$  maps to the plane $v = -c$ by $s \ra (16 (s+c)^2 s, (s+c)(c-3s), -c)$).

We conclude also easily that $ Sing (P) = \Gamma$. Since, if $p$ were another 
singular point of $P$, any plane through $p$ would intersect $P$ in a reducible
curve; but projection with centre $p$ yields a double cover of the plane
$\PP^2$, and a general line cannot be tangent to the branch locus,
 thus we get a contradiction.

Let now $X \subset \PP^3$ be a quartic surface which has a twisted cubic curve $\Gamma$ as
cuspidal curve: then it follows by the theorem \ref{edge} proved in the last section that
$X$ is unique, whence it coincides with the tangential developable of $\Gamma$. 

An alternative argument is as follows:
if two general polar surfaces 
$$ \Sigma _i y_i \partial X/ \partial x_i = 0,\  \Sigma _i z_i \partial X/ \partial x_i = 0$$
are shown to intersect along $\Gamma$ with multiplicity three, then the dual variety of $X$
is a curve $D$. 

Unfortunately, as pointed out by the referee, this statement is not obvious
 for a general $X$ as above, but indeed
for our explicit surface $P$ a direct calculation with Macaulay shows that the dual
variety of $P$ is a curve $D$ (and also that $D$ has degree $3$, 
but we do not need this fact).

Once we know that the dual variety of $X$ is a curve $D$, by biduality, 
$X$ is a developable surface which is not a cone,
so $X$ is a tangential developable, and its singular curve $\Gamma$ must be the edge of
regression. We conclude thus that $X$ is the tangential developable of $\Gamma$,
and that the dual variety $D$ of $X$ is the curve of osculating planes of $\Gamma$.
We conclude also that $D$ is then a twisted cubic curve, so $X$ is the dual surface
of a twisted cubic curve. $X$ is also projectively unique since $D$ is
projectively unique.

\QED

\section{The 2-dimensional picture}

In this section we shall assume that $u,v$ are local coordinates in the plane
(or local parameters for a two-dimensional local ring $A$), and that $a,b$ are 
local holomorphic functions at the origin, which  are invertible and take small values
 in a neighbourhood of the origin (respectively, $a,b$ are units of $A$).

Then we write $ b = c a$, and consider new coordinates $U,V$ such that $u = a^2 U$,
$ v = a^2 V$: in these new coordinates our discriminant $\Delta$ is divisible by $a^8$,
and after dividing by $- 4^4a^8$ we obtain the function 
$$ \delta_c:=  - U^2 V^2 - \frac{9}{8} UV (c)^2 + c^2 V^3 +  U^3 + \frac{27}{16^2}
 c^4  .$$
In other words, we could have assumed without loss of generality that $a=1$. 

We can further simplify the above equation by taking a cubic root
$\lambda$ of $c$, and considering new coordinates $u', v'$ with
$ U = \lambda^4 u_0 , V = \lambda^2 v_0$: then our equation, after dividing
by $\lambda^{12}$ becomes
$$ \delta:=  - u_0^2 v_0^2 - \frac{9}{8} u_0 v_0  +  v_0^3 +  u_0^3 + \frac{27}{16^2} .$$

Let us determine exactly the singular points of this curve, where for simplicity
of notation we replace $u_0, v_0$ by $u,v$ respectively. 
In other words, we could have assumed from the onset $a=b=1$, and we consider

$\delta(u,v) : = P (u,v,1,1) =  - u^2 v^2 - \frac{9}{8} u v  +  v^3 +  u^3 + 
\frac{27}{16^2},$

and our previous calculation of $u (\partial P / \partial u) - v (\partial P / \partial v) $
shows that the singular points satisfy $ u^3 = v^3$;
since however the origin does not lie on our curve, we may set, for such a 
singular point, $ u = \zeta v$, where $\zeta$ is a cubic root of $1$, and $ v \neq 0$.

We look now at $(\partial P / \partial u) = -2 \zeta v^3 
 - \frac{9}{8} v + 3 \zeta^2 v^2 = 0$, but disregarding the root $v=0$;
whence, we get  the equation $  v^2 
 + \frac{9}{16} \zeta^2  - \frac{3}{2} \zeta v = ( v -\frac{3}{4} \zeta)^2 = 0$.

Thus we conclude that the three cuspidal points are the three points
$$ u = \frac{3}{4} \zeta^2, v =\frac{3}{4} \zeta, (where \  \zeta^3 = 1). $$ 

Our goal is to understand the local braid monodromy of the curve 
$\delta = 0$, for a good projection given by a linear form $x$.

In order to understand the forthcoming calculations, observe that our curve 
$\Delta$ is the image of the ramification curve $ 4 zw = 1$, thus we have a degree
4 rational function $x$ on $\PP^1$, which must be branched on 6 points.

Three of these will correspond to the images of the 3 cusps, 2 will come from 
the two points at infinity where the line $z=0$ is tangent, whence there will be
exactly another branch point corresponding to a line $ x= x_0$ tangent at a smooth point.

Therefore the factors of the braid monodromy factorization will be four,
one half twist, and three cubes of a half twist. In the next section
we are going to calculate it for a suitable choice of coordinates.

\section{Making the  three cusps real.}

As we mentioned in the introduction, a three cuspidal quartic has as dual curve
a rational irreducible cubic 
(by Pl\"ucker's formulae its degree is $4 \times 3  - 3 \times 3 = 3$) which 
is by biduality nodal, since the dual of a cuspidal cubic is a cuspidal cubic.

Indeed the node is dual to the bitangent line at infinity. Our curve $\Delta$
will have two real tangents, and as a consequence three  flexes which are  not
all real, since
only one of the three cusps is real.

We easily construct however a quartic with three real cusps if we take the dual curve of
the affine curve
$$  D := \{ (X,Y)| F(X,Y)  = Y^2 - X^2 (X-1)  = 0\} $$
Since in homogeneous coordinates $(X,Y,Z)$ we have

$F(X,Y,Z)  = Y^2Z - X^2 (X-Z)  $, hence the gradient of $F$ is given by

$\nabla F = (-3 X^2 + 2 XZ, 2 YZ, X^2 + Y^2)$,

 in view of the standard 
parametrization of $D$ given by
$$ X = (t^2 + 1), Y = t (t^2 + 1) , Z=1  $$
(for which $ t = \infty$ goes to the point at infinity of $D$), we
get a parametrization of the dual curve $C$ as
$$( - 3 ( 1 + t^2) + 2, 2 t , (1 + t^2)^2 ) . $$   
The two flexes not at infinity occur for $ X = 4/3$, and $ t = \pm (3)^{-1/2}$,
and correspondingly we have the two flexes 
$(-2,  \pm 2 \ (3)^{-1/2} , (4/3)^2) = ( -\frac{9}{8}, \pm \frac{3}{8} \sqrt 3, 1)$
(the third flex occurs at the origin).

Using the above parametrization or using the computer we may calculate
the equation of $C$ as
$$ C := \{ (x,y) | (x^2 + y^2)^2 + x^3 + 9 x y^2 + \frac{27}{4} y^2 = 0 \}. $$
The advantage of this equation is that it is biquadratic in $y$, so $y^2$
is solution of the quadratic equation 
$$ (y^2)^2  + ( 2 x^2 + 9 x +  \frac{27}{4} ) y^2  + (x^3+ x^4) = 0$$
thus
$$ 2 y^2 = -   ( 2 x^2 + 9 x +  \frac{27}{4} ) \pm \sqrt 
{ 32 x^3 + 108 x^2 + \frac{27 \times 9}{2}
x + \frac{27^2}{16}  }, $$
and, if  we set $ A:=  ( 2 x^2 + 9 x +  \frac{27}{4} )$ ,
the  discriminant  $\Theta = A^2 - 4 (x^3+ x^4)$ appearing in the above square root 
 is clearly positive
for $ x >0$, and clearly vanishes for $ x = -\frac{9}{8}$.
On the other side, twice the derivative of $\Theta$ equals  
$ 192 x^2 + 432 x + 243 = 3 (8x + 9)^2 \geq 0$, 
thus $\Theta$ is strictly monotone, whence $\Theta$ is positive exactly
for $ x > -\frac{9}{8}$.

To complete the picture, observe that $\Theta > A^2$ iff $x^4 + x^3 < 0 $,
i.e., iff $ -1 < x < 0 $, and that 
$ A(x) = 2 ( (x+ \frac{9}{4})^2 - \frac{27}{16})$ is positive exactly
outside the interval $ \frac{3}{4} (- 3  \pm \sqrt3)$, and 
note that $\frac{3}{4} (3 - \sqrt3) < 1.$

We have thus the following picture:

\font\mate=cmmi8
\def\min{\hbox{\mate \char60}\hskip1pt}

\begin{center}
\begin{picture}(180,220)(0,-40)

\put(0,80){\line(1,0){140}}

 \put(110,0){\line(0,1){160}}

\put(20,80){\circle*{3}}

\put(30,80){\circle*{3}}

\put(20,131){\circle*{3}}

\put(20,29){\circle*{3}}

\put(35,87){-1}

\put(10,87){-$\frac{9}{8}$}

 \put(115,87){0}

\put(110,80){\circle*{3}}

\put(160,80){x-axis}

\put(100,170){y-axis}


%
\qbezier[40](20,131)(40,80)(20,29)

\qbezier[40](20,131)(40,80)(110,80)

\qbezier[40](20,29)(40,80)(110,80)

%
\end{picture}

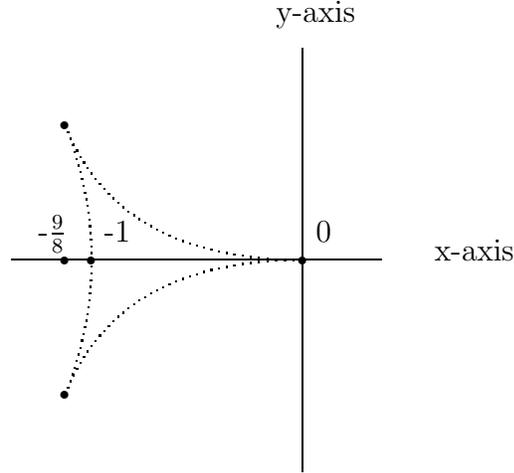
The real part of the curve C.
\end{center}

Using the previous description, and looking at the above picture, we may now 
easily describe the motion of the roots $y$ as $x$ moves along the real axis
from $ + \infty$ to $-\frac{9}{8} $. 

For $x> 0$ we have exactly 4 imaginary roots 
$ A_1(x) =  i Y_1(x), A_2(x) = i Y_2(x), B_2(x)=  -  i Y_2(x), 
B_1(x) = - i Y_1(x)$, where $ Y_1(x), Y_2(x) \in \R$ and $Y_1(x) > Y_2 (x) > 0$.

For $x=0$, $A_2, B_2$ become $ = 0$, and then for $ -1 < x < 0 $ $A_2, B_2$
 become real and opposite, $ B_2$ is positive and grows as $x$ decreases, 
while $A_1$ remains imaginary and with decreasing absolute value.

For $ x = - 1$ $A_1, B_1$ become $ = 0$, while $ B_2 =\frac{1}{4}, 
 A_2 = - \frac{1}{4}$, finally in the interval $ - \frac{9}{8} < x < -1$
we have 4 distinct real roots $ B_2(x) > B_1(x) > 0 >   A_1(x)  = -  B_1(x)
>  A_2(x)  = -  B_2(x)$ and both roots $ B_2(x) > B_1(x)$ grow as $x$ approaches
$ - \frac{9}{8}$: for this value we have $ B_2 = B_1 =\frac{3}{8} \sqrt 3.$

\section{Braid monodromy of $C$ and fundamental group of the complement}

We take as projection of the pair $ (\C^2, C)$ a linear form very close to $x$,
since the projection $x$ has a double critical value $ x= - \frac{9}{8} $:
the effect is then that we split into two factors the corresponding braid.

In order to get a simple picture let us take as base point
$x_0=\frac{3}{4}(\sqrt{3}-3)$. Then $-1<x_0<0$, and we have two real
roots $A_2$, $B_2$ and two purely imaginary roots $A_1$, $B_1$ as
in Figure 2, left side.

  If we
move on the real axis to the right to $x=0$ then the
two real roots meet at zero and we
have a  horizontal vanishing arc connecting the real roots with a
$3/2$ -twist around this arc as local monodromy. If we move on the real
axis to the left
to $x= -1$ then the two complex roots meet at zero and we get a
vertical vanishing arc connecting the imaginary roots with one
half-twist around this arc as local monodromy. If instead of going
directly to $x=-1 $ we make a half turn around it counterclockwise and continue
  along the real axis to the left then the imaginary roots turn
  $  \pi /4$ counterclockwise around zero,
  the top root $A_1$ becomes real negative and "runs after $A_2$"
  and the bottom root $B_1$ becomes real positive and "runs after
  $B_2$". The roots meet
  for the appropriate critical
values of x (in a neighbourhood of $ - \frac{9}{8}$).

  The corresponding vanishing arcs
are just straight intervals connecting $A_1$ with $A_2$ and $B_1$
with $B_2$. They are disjoint so the order of these two critical
values of x is not important. The corresponding monodromy factors
are 3/2-twists around these arcs and they commute.We take a cyclic
order of the paths counterclockwise around $x_0$ starting first with the
two critical values near $-\frac{9}{8}$, then proceeding with
  the critical value '$x=-1$' (actually , near  $x=-1$: recall in fact
that we changed
slightly the axis of projection to split the two critical values of $x$ near
$-\frac{9}{8}$, and finally ending with '$x=0$'.

The base point is usually chosen far away from the critical
values. We should do this also in our case since we deal with a local
picture and we may at a later convenience want  to relate it to a 
global picture
where other
monodromy factors occur, coming from other critical values.

  We  move then the base
point to the right along the real axis, passing $x=0$ on the
right, making a half-turn clockwise around it. The paths from $x_0$
to the critical values of x will be dragged along. The real roots
$A_2$ and $B_2$ get closer to zero and then turn clockwise around
zero by a $(3/2) \pi$ turn (in fact, a full turn of $x$ around $x=0$ produces
a $3/2$-turn
of the roots). The whole picture of vanishing arcs moves with
them. We get a configuration as on the right side of Figure 2.

When $x$ moves further to the right along the real axis the roots
remain purely imaginary and just move further away from zero.

We have thus obtained a complete determination of the braid
monodromy of $C$, which is illustrated by Figure 2 and
summarized in the following

\begin{teo}       Consider the arcs depicted on the right side
of Figure 2.

Then the braid monodromy of $C$ is given,  in this order, by the cube of a
half-twist around the arc connecting $A_1$ and $A_2$, by the cube of
a half-twist around the arc connecting $B_2$ and $B_1$, by the
half-twist around the arc connecting $A_1$ and $B_1$ and finally by the cube
of a half-twist around the arc connecting $A_2$ and $B_2$.

\end{teo}

   \begin{center} \begin{picture}(300,160)

\put(30,80){\line(1,0){80}}

\put(30,80){\line(1,1){40}}

\put(70,40){\line(0,1){80}}

\put(70,40){\line(1,1){40}}

\put(220,74){\line(0,1){12}}

\put(30,80){\circle*{3}}

\put(110,80){\circle*{3}}

\put(70,40){\circle*{3}}

\put(70,120){\circle*{3}}

\put(65,130){$A_1$}

\put(65,25){$B_1$}

\put(10,80){$A_2$}

\put(120,80){$B_2$}

\put(220,50){\circle*{3}}

\put(220,74){\circle*{3}}

\put(220,86){\circle*{3}}

\put(220,110){\circle*{3}}

\put(220,85){\oval(50,50)[r]}

  \put(220,75){\oval(50,50)[l]}

  \put(220,92){\oval(36,36)[r]}

  \put(220,68){\oval(36,36)[l]}

  \put(220,90){\oval(20,20)[r]}

  \put(220,70){\oval(20,20)[l]}

\put(215,120){$A_1$}

\put(215,35){$B_1$}

\put(217,63){$A_2$}

\put(203,88){$B_2$}

\put(220,74){\line(0,1){12}}

\end{picture} \\

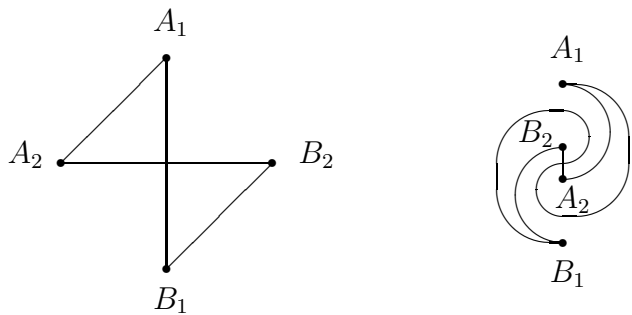
   Figure 2. {\em Vanishing arcs}

\end{center}

Let us now
consider again the left part of Figure 2 and let us take as base
point for the fundamental group of the fibre
  $\C - \{ A_1, A_2, B_2, B_1 \}$ a point
$y_0$ with large positive real part and small positive imaginary part.

We consider then a geometric basis $\alpha_1, \alpha_2, \beta_2,\beta_1$
of $ \pi_1 (\C - \{ A_1, A_2, B_2, B_1 \}, y_0)$,
where  for instance $\alpha_1$ is given
  by  a subsegment $s_{A_1}$ on the  segment joining $y_0$ with $A_1$,
  followed by a full small circle around $A_1$, and then followed by the
inverse path of $s_{A_1}$ (the other loops are defined similarly).

By the van Kampen theorem the fundamental group $ \pi_1 (\C^2  - C, y_0)$
is generated by $\alpha_1, \alpha_2, \beta_2,\beta_1$ subject to the relations
coming from the braid monodromy, thus we obtain a presentation
$$ \pi_1 (\C^2  - C, y_0) = < \alpha_1, \alpha_2, \beta_2,\beta_1|
\alpha_1 \alpha_2 \alpha_1 = \alpha_2 \alpha_1 \alpha_2,
\beta_1 \beta_2 \beta_1 = \beta_2 \beta_1 \beta_2,$$

$$\alpha_2 \beta_2 \alpha_2  = \beta_2 \alpha_2 \beta_2,
 \beta_2 \beta_1 = \alpha_1 \beta_2 >.$$

We need only explain the last relation, coming from the relations
$\sigma(\gamma) = \gamma$, where $\sigma$ is the half twist on the curve
$\tau$, corresponding
to the  vertical tangency for $x= -1$, which makes the two roots $A_1, B_1$ become equal.

The action of this half twist $\sigma$ is the following
\begin{itemize}
\item
$ \alpha_1 \ra \alpha_1  \beta_2 \beta_1  \beta_2^{-1} \alpha_1^{-1}$
\item
$ \alpha_2 \ra \alpha_1  \beta_2 \beta_1^{-1}  \beta_2^{-1} \alpha_2 
\beta_2 \beta_1 \beta_2^{-1} \alpha_1^{-1}$
\item
$\beta_2 \ra \beta_2$
\item
$\beta_1 \ra \beta_2^{-1} \alpha_1 \beta_2.$

\end{itemize}

We end by describing the monodromy homomorphism $\mu$ of the degree 4 covering 
$ \pi_1 (\C^2  - C, y_0) \ra \SSS_4$.
The action of $ \mu$ must then be, up to conjugation
in $\SSS_4$, the following one:

\begin{itemize}
\item
$ \alpha_1 \ra (1,2)$
\item
$ \alpha_2 \ra (2,3)$
\item
$\beta_2 \ra (2,4)$
\item
$\beta_1 \ra (1,4).$

\end{itemize}

In fact, we observe first that

(**) the generators $\alpha_1 ,\alpha_2,\beta_2 ,\beta_1 $  must map to transpositions,
 and the image of $\mu$ is transitive.

One sees then, using the presentation of $\pi_1(C^2 - C)$, that 

(***) there is a unique 
homomorphism of $\pi_1(C^2 - C)$ into $\SSS_4$ (up to
conjugation in $\SSS_4$) which satisfies (**).

If instead we want to calculate the fundamental group of the complement
$ \pi_1 (\PP^2  - C), y_0)$ we must add the relation 
$ \alpha_1 \alpha_2 \beta_2 \beta_1 = 1$ and we can then
 simplify the presentation obtaining 
$ \pi_1 (\PP^2  - C, y_0) =< \alpha_1, \alpha_2 |  
\alpha_1 \alpha_2 \alpha_1 = \alpha_2 \alpha_1 \alpha_2,\ 
\alpha_2 \alpha_1  \alpha_1 \alpha_2 = 1>$, which is the 
spherical braid group of three points in $\PP^1$, as shown
long ago by Zariski, and also in greater 
generality by Moishezon (cf. \cite{zar}, \cite{moi1}).

This calculation shows that the degree four covering is also branched
on the line at infinity, where the local monodromy is the double
transposition $(1,3)(2,4)$ (in the Galois case we have three branch lines,
corresponding to the three nontrivial elements of $(\Z/2)^2$,
and we have the standard model for a bidouble cover given by a special
projection of the Veronese surface, as described in \cite{cat99}, page 100.

Let us briefly recall it: consider the Veronese surface $V$ , i.e., the
variety of symmetric matrices of rank 2
$$rank 
\begin{pmatrix}
  x_1 & w_3& w_2 \\
  w_3 & x_2 &w_1  \\
  w_2 & w_1 & x_3  \\

\end{pmatrix} = 1.$$
$V$ is isomorphic to $\PP^2$ with coordinates $(y_1, y_2, y_3)$ by setting
$$ x_i = y_i^2, \ \  w_1 = y_2 y_3, w_2 = y_1 y_3, w_3 = y_1 y_2,  $$
and the projection $\pi :V \ra \PP^2$ given by $(x_1, x_2, x_3)$
corresponds to the $(\Z/2)^2$ Galois cover
$$(y_1, y_2, y_3) \ra  (y_1^2, y_2^2, y_3^2). $$

The deformed degree 4 covering is then 
a  slightly less special, yet very interesting projection of the Veronese surface.

\section{The tangential developable $F$ of the twisted cubic $\Gamma$.}

This final section is devoted to the proof of an interesting characterization 
of the above surface

\begin{teo}\label{edge}
The tangential developable $F$ of the twisted cubic $\Gamma$ is the unique irreducible
surface of degree $4$ in $\PP^3$ which has the twisted cubic $\Gamma$ as a cuspidal curve.
\end{teo}

\begin{oss}
1) After proving the theorem, we looked again with more care at
the book
by Conforto and Enriques " Le superficie razionali" (\cite{supraz}),
 where Chapter V  is mostly devoted to
the surfaces of degree $4$ which contain $\Gamma$ as a double curve.
On page 114 the above theorem is mentioned, and a different proof is briefly sketched
in a footnote.

It is then mentioned in \cite{supraz} that the complete classification
of surfaces of degree $4$ ruled by lines, started by Calyley in \cite{cay}.
 was achieved by Cremona (\cite{crem}),
and a later classification was also given by G.Gherardelli (\cite{gher}).
For lack of time (pending deadline),
 we are not in a position to determine who gave the first
proof of the above theorem \ref{edge}.

We hope however that our modern description of quartic surfaces having the
twisted cubic as double curve may be found simple and useful.

2) Note that, if we view $\PP^3$ as the space of effective degree $3$ divisors on $\PP^1$,
then $\Gamma$, $ F - \Gamma$, $\PP^3 - F$ are exactly the three $\PP GL(2)$ -orbits,
corresponding to the divisors of respective types $ 3 P$, $ 2 P_1 + P_2$ ($ P_1 \neq P_2$),
 $ P_1 + P_2 + P_3$ with $ P_1 , P_2 , P_3$ three distinct points.
\end{oss}

\Proof {\em of Theorem \ref{edge} }
The twisted cubic curve $\Gamma$ is the image of $\PP^1$ under the third Veronese mapping
$ v_3 (t_0,t_1) : = ( t_0^3,t_0^2 t_1, t_0t_1^2, t_1^3)$, and its projective coordinate 
ideal $I_{ \Gamma} $ is generated by three quadrics, the determinant of the 
$2 \times 2$-minors of the matrix
$$ A = \begin{pmatrix}
  x_0 & x_1 & x_2  \\
  x_1 & x_2 & x_3  \\

\end{pmatrix} .$$
Thus we have a three dimensional vector space $V$ consisting of the quadrics containing 
 $\Gamma$. Thus $ V : = I_{ \Gamma}(2) = H^0 (\PP^3, \I_{ \Gamma}(2))$, where
$\I_{ \Gamma}$ is the ideal sheaf of $\Gamma$, is generated by
$$ Q_0 : =  x_1 x_3 -  x_2^2 , \  Q_1 : =  - x_0 x_3  +   x_1  x_2, 
\  Q_2 : =  x_0 x_2 -x_1^2 .$$ 
Observe now that, to each point  $x \in \Gamma$ is associated a unique quadric cone
$Q_x $ containing $ \Gamma$ and with vertex $x$: since projection with centre $x$
maps $\Gamma$ to a plane conic.
Thus we have a map $\Gamma \ra \PP (V)$ associating $Q_x$ to $x$ (in our notation,
$\PP (V)$ denotes the set of $1$-dimensional subspaces of the vector space $V$).
We denote by $\tilde{\Gamma}$ the image of this map, and observe that there
is thus a canonical bijection between ${\Gamma}$  and $\tilde{\Gamma}$ .

\begin{lem}
Consider in the projective plane $\PP (V) : = \{ Q \ | \Gamma \subset Q\}$
the quartic curve $ \{ Q | det Q = 0\}$. Then this quartic curve 
is the conic $\tilde{\Gamma} \subset \PP (V)$ counted with multiplicity $2$.
\end{lem}

\Proof {\em of the Lemma. }
Observe that if a quadric $Q$ contains $\Gamma$, then $ Rank (Q) \geq 3$ since
$\Gamma$ is irreducible.

{\bf CLAIM:} If  $ Rank (Q) = 3$ (then $Q$ is a quadric cone)  the vertex $x$
of $Q$ is a point of $\Gamma$. 

{\bf Proof of the claim}: otherwise $\Gamma$ would be a Cartier divisor on $Q$,
whence it is known (but we reprove it below) that its degree should  be even.

Let in fact $\F_2$ be 
the blow-up of $Q$ at $x$, thus a basis of $ Pic (\F_2)$ is given by the excetional
curve $\sigma$, and by the strict transform $F$ of a line, which satisfy
$ F^2=0, \ \sigma^2 = -2 ,  \ \sigma F = 1$. The plane section $H$ is linearly
equivalent to $ 2 F + \sigma$, and let $\Gamma \equiv a \sigma + b F$.
From the equations $ \sigma \Gamma = 0, H \Gamma = 3$ we obtain $ b - 2a =0,b = 3$,
a contradiction (in reality the class of $\Gamma$ is  $ 3 F + \sigma$).
\qed

It follows that the quartic curve $ \{ Q \in \PP (V) | det Q = 0\}$ is 
set theoretically the curve $\tilde{\Gamma}$, which is a rational curve.
But $\tilde{\Gamma}$ is homogeneous by the action of $\PP GL(2)$ acting on
$\PP (V)$, thus $\tilde{\Gamma}$ is smooth and must be a conic: the assertion follows then
right away.

\qed \ for the Lemma

We analyse now the space of quartics which have $\Gamma$ as a double curve:

\begin{lem}
$ U : = H^0 (\PP^3, \I_{ \Gamma}^2(4))  \cong Sym^2 (V) .$
\end{lem}
\Proof {\em of the Lemma. }
Let $F \in H^0 (\PP^3, \I_{ \Gamma}^2(4))$ : we need to show that $F$ is equal to a quadratic 
polynomial $f (Q_0,Q_1,Q_2)$. Let us first consider the divisor cut by $F$ on the smooth
quadric $Q_1$: since $ Pic (Q_1)$ has as basis the respective rulings $L_1$ and $L_2$,
 and the hyperplane divisor $H$ is linearly equivalent
to $ L_1 + L_2$, we see by direct calculation that, under the isomorphism
$Q_1 \cong \PP^1 \times \PP^1$,  $ div_{Q_1} (Q_0) = \Gamma + div (u_0)$,
 $ div_{Q_1} (Q_2) = \Gamma + div (u_1)$, where $(u_0, u_1)(v_0, v_1)$ are suitable
 coordinates
on $\PP^1 \times \PP^1$. 

It follows in particular that the quadric cones
$Q_x$ cut on $Q_1$ the curve $\Gamma$ plus the line of the first ruling passing through $x$,
and more importantly that $  div_{Q_1} (F) = 2 \Gamma + div ( \phi(u_0, u_1))$,
where $\phi$ is a quadratic polynomial.
Whence, $  div_{Q_1} (F) =  div (\phi(Q_0, Q_2))$, and there is a quadratic form $Q_3$
in $\PP^3$ such that
$$ F = \phi(Q_0, Q_2) + Q_1 Q_3.$$
Since however $\phi(Q_0, Q_2) \in H^0 (\PP^3, \I_{ \Gamma}^2(4))$, it follows
that $Q_1 Q_3\in H^0 (\PP^3, \I_{ \Gamma}^2(4))$ and thus $Q_3 
\in H^0 (\PP^3, \I_{ \Gamma}(2))$, so that $Q_3$ is a linear combination
of $ Q_0,   Q_1 , Q_2$.

\qed \ for the Lemma

For such a surface $F$ as above, $\Gamma$ is a double curve, and we are going to show 
that a general such surface possesses $4$ pinch points on $\Gamma$.
As a first step in this direction, we describe the conormal bundle to $\Gamma$.

\begin{lem}
$ \NNN^*_{\Gamma}(2) \cong \hol_{\PP^1}(1) \oplus  \hol_{\PP^1}(1) .$
\end{lem}
\Proof {\em of the Lemma. }

We know that $ Q_0,   Q_1 , Q_2 \in H^0 (\PP^3, \I_{ \Gamma}(2))$ induce,
under the surjection $\I_{ \Gamma}(2) \ra \I_{ \Gamma}/\I_{ \Gamma}^2 (2)
= \NNN^*_{\Gamma}(2) $,
three sections $ q_0,   q_1 , q_2$ which generate the Rank $2$ bundle
$\NNN^*_{\Gamma}(2)$. Since $ det (\NNN^*_{\Gamma}(2))$ has degree $2$
on $\PP^1$, as it is easily seen by the cotangent bundle sequence
for $ \Gamma \subset \PP^3$, it follows that $\NNN^*_{\Gamma}(2)$ splits
either as $\hol_{\PP^1}(1) \oplus  \hol_{\PP^1}(1) $ or 
$\hol_{\PP^1}(0) \oplus  \hol_{\PP^1}(2) $.

But we can exclude the second case since any quadric cone $ Q_x\in V$
induces a section $q \in H^0 (\NNN^*_{\Gamma}(2) )$ whose two components
have a simple zero at $x$, and no other zero. 

\qed \ for the Lemma

We have now a discriminant map 
$$ \delta : Sym^2 ( \NNN^*_{\Gamma}(2) ) =
\hol_{\PP^1}(2) \oplus  \hol_{\PP^1}(2) \oplus  \hol_{\PP^1}(2) 
\ra \hol_{\PP^1}(4) ,$$
given by $ \delta (a,b,c) = 4 ac - b^2  $ and vanishing on the simple tensors.

Assume now that $ q_j = (a_j, b_j)$: then we get an equation for the pinch points
of a surface
$$ F = \sum_{i,j} \lambda _{i,j} Q_i Q_j,$$
namely,
$$ \Delta (F) = 4 ( \sum_{i,j} \lambda _{i,j} a_i a_j) 
( \sum_{i,j} \lambda _{i,j} b_i b_j) - ( \sum_{i,j} \lambda _{i,j} 
[a_i b_j + b_i a_j ])^2 \in H^0 (\PP^1, \hol_{\PP^1}(4)).$$ 

We see then that a general surface $ F \in \PP (  Sym^2 (V))$ has $4$ pinch points
along $\Gamma$, and that 

{ \bf (***) $\Gamma$ is a cuspidal curve for $F$ if and only if $\Delta(F) \equiv 0$.}

We observe then that (***) is a system of $5$ quadratic equations vanishing
on the Veronese surface ( the surfaces in $\PP (  Sym^2 (V))$ which are squares
$Q^2$ of $ Q \in V$). We can then bet that our problem is equivalent to
the problem of the number of conics tangent to $5$ fixed lines in the plane.

We show that this is indeed the case, because first of all

1) $\PP (  Sym^2 (V))$ is the space of conics in $\PP ( V^{\vee})$.

2) $\Gamma$ is a cuspidal curve for $F$ if and only if $\Delta(F)$ vanishes in $5$
fixed distinct points $x(1), \dots , x(5) \in \Gamma$, i.e., 
$x(1), \dots , x(5) \in \Gamma$ are pinch points for $F$.

3) the following holds:

\begin{lem} $x \in \Gamma$ is a pinch point for $F$ if and only if the corresponding
conic $C_F$ in $\PP ( V^{\vee})$ is tangent to the line in $\PP ( V^{\vee})$
dual to the point $\tilde{x} \in \tilde{\Gamma}$ corresponding to $x$ 
(i.e., $\tilde{x} $ is the point
given by the quadric cone $Q_x$). 

\end{lem}
\Proof {\em of the Lemma. }

Let us take coordinates $(t_0, t_1)$ in $\PP^1$ such that $x$ 
corresponds to the point $t_1=0$. Likewise, since $\tilde{\Gamma}$ is a smooth conic,
we may take a basis of $V$ corresponding to quadric cones $ Q_0,   Q_1 , Q_2$
with vertices in the respective points of $\Gamma$ corresponding to
$ t_0 = 0, t_1 =0, t_2 : = t_0 - t_1 =0$.

Observe then for later use  that the point $x$ corresponds to the point 
$\tilde{x}  = (0,1,0)$ in 
$\PP (V)$, so its dual line, in the dual basis coordinates, will be the line
$ y_1 = 0$.

The evaluation of the three sections $ q_0,   q_1 , q_2$ in the fibre of
the conormal bundle at $x :  t_1=0, t_0 = 1$ yields respective 
vectors forming a matrix

$$ A'' : = \begin{pmatrix}
  1 & 0 & 1  \\
  b_0 & 0 & b_2  \\

\end{pmatrix} .$$
Then a quartic $ F = \sum_{i,j} \lambda _{i,j} Q_i Q_j$ has our point as a pinch
point if and only if the following $ 2 \times 2 $ symmetric matrix has
zero determinant:

$$ M : = \begin{pmatrix}
  \lambda_{00} + \lambda_{02} + \lambda_{22}  & \lambda_{00}b_0 +  \frac{1}{2} \lambda_{02}
(b_0 + b_2) + \lambda_{22} b_2  \\
  \lambda_{00}b_0 +  \frac{1}{2} \lambda_{02}
(b_0 + b_2) + \lambda_{22} b_2 & \lambda_{00}b_0^2 + \lambda_{02}b_0 b_2 + \lambda_{22} b_2^2  \\

\end{pmatrix} .$$

Then we have 

$$ det (M) = (b_0 - b_2)^2 ( \lambda_{00} \lambda_{22} - 4 \lambda_{00}^2).$$
We then observe that $(b_0 - b_2)^2  \neq 0$ because not all the above
quartics $F$ have a pinch point in $x$. 

Therefore, the condition that $x$ be a pinch point is exactly given by
$( \lambda_{00} \lambda_{22} - 4 \lambda_{00}^2) = 0$, i.e., by the condition that
the conic $ C_F = \{(y_0,y_1,y_2) | \sum_{i,j} \lambda _{i,j} y_i y_j = 0 \}$
be tangent to the line $\{ y_1=0 \}$ dual to the point $\tilde{x}$.

\qed \ for the Lemma

We are now ready to finish the proof of  Theorem 6.1: assume that $F$ is a quartic
surface which has $\Gamma$ as a double curve, and assume that $F$ is not the square of
a quadric $Q \in V$. Then the corresponding conic $C_F$ has rank $\geq 2$.

Assume further that $\Gamma$ is a cuspidal curve for $F$: this holds if and only if the 
conic $C_F$ is tangent to five fixed lines $L_1, \dots , L_5$ dual to $5$ points
$ \tilde{x}_1, \dots , \tilde{x}_5 \in \tilde{\Gamma}$.

Then the dual  conic $C_F^{\vee}$ passes through the five points 
$ \tilde{x}_1, \dots , \tilde{x}_5 \in \tilde{\Gamma}$, and therefore coincides with
$\tilde{\Gamma}$.

We have thus shown that there is one and only one such quartic $F$ which is not
a quadric counted with multiplicity $2$, so we conclude that $F$, which is a union of
$\PP GL(2)$ orbits, is exactly the tangential developable surface of
$\Gamma$, which is an irreducible surface.

\qed

\begin{footnotesize}
\noindent

{\bf Note.}
This article is essentially based on classical  mathematics, and it is accordingly 
written in
classical style, sometimes referred to by referees as: 
"a style which is suitable for Conference Proceedings".
We hope that the style may be suitable for the reader.

{\bf Ackowledgements.}

We would like to thank the referee for pointing out a gap in the original proof of 
Theorem \ref{disc}, and correcting a couple of minor mistakes.
We also thank Fabio Tonoli for providing a Macaulay 2 Script verifying that the
dual variety of $P$ is indeed a curve.

Finally, the first named author would like to acknowledge the hospitality 
of M.S.R.I. in march 2004, where the first calculations of the braid monodromy
were begun. 
\end{footnotesize}

{\bf Dedication.}

Finally, the  manifold occurrence of words such as
"Veronese surface" or "Veronese embedding" points out
 the appropriateness of this article to celebrate
the 150-th anniversary of the birth of Giuseppe Veronese.

\vfill

\bigskip

\noindent 
Prof. Fabrizio Catanese\\
Lehrstuhl Mathematik VIII\\
Universit\"at Bayreuth, NWII\\
 D-95440 Bayreuth, Germany

e-mail: Fabrizio.Catanese@uni-bayreuth.de

\noindent
Prof. Bronislaw Wajnryb \\
Department of Mathematics \\
Technion\\
  32000 Haifa, Israel

e-mail: wajnryb@techunix.technion.ac.il

\end{document}